\documentclass{article}
\usepackage{amsmath, amsfonts, amssymb, amsthm, a4wide}
\title{A Recent Look at the Quantum Berezinian in the Yangian $Y(\mathfrak{gl}_{m|n})$}
\author{Lucy Gow\footnote{lucyg@maths.usyd.edu.au} \\ \normalsize{University of Sydney}}
\begin{document}
\maketitle
\newcommand{\pa}[1]{\overline{#1}}
\newtheorem{lemma}{Lemma}[section]
\newtheorem{proposition}[lemma]{Proposition}
\newtheorem{defn}{Definition}[section]
\newtheorem{theorem}{Theorem}
\newtheorem*{dualform}{Dual Form of the Quantum Berezinian}
\newtheorem{remark}{Remark}[section]
Brundan and Kleshchev recently introduced a new family of presentations of the Yangian $Y(\mathfrak{gl}_{n})$ associated to
 the general linear Lie algebra $\mathfrak{gl}_n$, and thus provided a fresh approach to its study. 
In this article, we would like to show how some of their ideas can be fruitfully extended to consider the Yangian 
$Y (\mathfrak{gl}_{m|n})$ associated to the Lie superalgebra $\mathfrak{gl}_{m|n}$.  In particular, we give a new proof of the result by Nazarov
that the quantum Berezinian is central.\footnote{The author would like to acknowledge the guidance and support of her PhD supervisor Alexander Molev.}
\section{Definition of $Y (\mathfrak{gl}_{m|n})$}
The Yangian $Y (\mathfrak{gl}_{m|n})$ is defined in \cite{Nazarov} to be the $\mathbb{Z}_2$-graded associative algebra over $\mathbb{C}$ with generators $t_{ij}^{(r)}$ and certain relations
described below.  We define the formal
power series
\begin{equation*}
t_{ij} (u) = \delta_{ij} + t_{ij}^{(1)} u^{-1} + t_{ij}^{(2)}u^{-2} + \ldots,
\end{equation*}
and a matrix
\begin{equation}\label{T}
T(u) = \sum_{i,j =1}^{m+n} t_{ij} (u) \otimes E_{ij} \; (-1)^{\pa{j} (\pa{i} +1)},
\end{equation}
where $E_{ij}$ is the standard elementary matrix and $\pa{i}$ is the parity of the index $i$. 
In analogy with the usual Yangian $Y(\mathfrak{gl}_{n})$ (see for example \cite{CP}, \cite{MNO}, \cite{Molev}), the defining relations are then expressed by the matrix product
\begin{equation*}
R(u-v) T_1 (u) T_2 (v) = T_2 (v) T_1 (u) R(u-v)
\end{equation*}
where
\begin{equation*}
R(u-v) = 1 - \frac{1}{(u-v)} P_{12}
\end{equation*}and $P_{12}$ is the permutation matrix:
\begin{equation*}
P_{12} = \sum_{i,j =1}^{m+n} E_{ij} \otimes E_{ji} (-1)^{\pa{j}}.
\end{equation*}
Then we have the following equivalent form of the defining relations:
\begin{equation*}
[t_{ij}(u),\;t_{kl}(v)] \;=\; \frac{(-1)^{\pa{i}\pa{j} +\pa{i} \pa{k} + \pa{j}\pa{k}}}{(u-v)} ( t_{kj}(u) t_{il} (v) - t_{kj} (v) t_{il} (u) ).
\end{equation*}
Throughout this article we will observe the following notation for entries of the inverse matrix of $T(u)$: 
\begin{equation*}
T(u)^{-1} =: \left(t'_{ij}(u) \right)_{i,j=1}^{n}.
\end{equation*}
A straightforward calculation yields the following relation in $Y(\mathfrak{gl}_{m|n})$:
{\small
\begin{eqnarray}\label{useful}
[t_{ij}(u), t'_{kl}(v)]
= \frac{(-1)^{\pa{i}\pa{j} +\pa{i}\pa{k} + \pa{j}\pa{k}}}{(u-v)}\cdot (\; \delta_{kj} \sum_{s=1}^{m+n} t_{is}(u) t'_{sl} (v)
-\;\delta_{il} \sum_{s=1}^{m+n} t'_{ks}(v)t_{sj}(u))  .
\end{eqnarray}
}
\section{Gauss Decomposition of $T(u)$}
In \cite{BK}, the Drinfeld presentation is described in terms of the quasideterminants of Gelfand and Retakh (\cite{GGRW}, \cite{GR}).  In this article we make use of the
analogous set of generators of the Yangian $Y(\mathfrak{gl}_{m|n})$.
First we recall the definition of the quasideterminants and some conventional notation.
\begin{defn}
Let $X$ be a square matrix over a ring with identity such that its inverse matrix $X^{-1}$ exists, and such that its $ji$th entry is an invertible element of the ring.  Then the $ij$th
\emph{quasideterminant} of $X$ is defined by the formula
\begin{equation*}
|X|_{ij} = \left((X^{-1})_{ji}\right)^{-1}.
\end{equation*}
\end{defn}
It is sometimes convenient to adopt the following alternative notation for the quasideterminants: 
\begin{equation*}
|X|_{ij} =: \left| \begin{array}{ccccc} x_{11} & \cdots & x_{1j} & \cdots & x_{1n}\\
&\cdots & & \cdots&\\
x_{i1} &\cdots &\boxed{x_{ij}} & \cdots & x_{in}\\
& \cdots& &\cdots & \\
x_{n1} & \cdots & x_{nj}& \cdots & x_{nn}
\end{array} \right|.
\end{equation*}
The matrix $T(u)$ defined in (\ref{T}) has the following Gauss decomposition in terms of quasideterminants (by Theorem 4.96 in \cite{GGRW}; see \S 5 in \cite{BK}):
\begin{equation*}
T(u) = F(u) D(u) E(u)
\end{equation*}
for unique matrices 
{\small
\begin{eqnarray*}
D(u) &=& \left( \begin{array}{cccc} d_1 (u) & &\cdots & 0\\
& d_2 (u) &  &\vdots\ \\
\vdots & &\ddots &\\
0 &\cdots &  &d_{m+n} (u)
\end{array} \right),\\
E(u) &=&  \left( \begin{array}{cccc} 1 &e_{12}(u) &\cdots & e_{1,m+n}(u)\\
&\ddots & & e_{2, m+n}(u)\\
& &\ddots & \vdots\\
0 & & &1 
\end{array} \right),\\
F(u) &=& \left( \begin{array}{cccc} 1 & &\cdots &0\\
f_{21}(u) &\ddots & &\vdots\\
\vdots & & \ddots & \\
f_{m+n,1}(u) & f_{m+n, 2}(u) &\cdots &1
\end{array} \right),
\end{eqnarray*}
}
where
{\small
\begin{eqnarray*}
d_i (u) &=& \left| \begin{array}{cccc} t_{11}(u) &\cdots &t_{1,i-1}(u) &t_{1i}(u) \\
\vdots &\ddots & &\vdots \\
t_{i1}(u) &\cdots &t_{i,i-1}(u) &\boxed{t_{ii}(u)}
\end{array} \right|, \\
e_{ij}(u) &=& d_i (u)^{-1} \left| \begin{array}{cccc} t_{11}(u) &\cdots &t_{1,i-1}(u) & t_{1j}(u) \\
\vdots &\ddots &\vdots & \vdots \\
t_{i-1,i}(u) &\cdots &t_{i-1,i-1}(u) & t_{i-1,j}(u)\\
t_{i1}(u) &\cdots &t_{i,i-1}(u) &\boxed{t_{ij}(u)}
\end{array} \right|,\\
f_{ji}(u) &=& \left| \begin{array}{cccc} t_{11}(u) &\cdots &t_{1, i-1}(u) & t_{1i}(u) \\
\vdots &\ddots &\vdots &\vdots \\
t_{i-1,1}(u) &\cdots &t_{i-1,i-1}(u) &t_{i-1,i}(u)\\
t_{ji}(u) &\cdots &t_{j, i-1}(u) &\boxed{t_{ji}(u)} 
\end{array} \right|
d_{i}(u)^{-1}.
\end{eqnarray*}
}It is easy to recover each generating series $t_{ij}(u)$ by multiplying together and taking commutators of $d_i (u); \; 1 \le i \le m+n$, and 
$e_i (u):= e_{i,i+1} (u), \,f_{i}(u)=f_{i+1,i}(u);~1 \le i < m+n$ (see \S 5 of \cite{BK}).
Thus the Yangian $Y(\mathfrak{gl}_{m|n})$ is generated by the coefficients of the latter.

\subsection{Some Useful Maps}
Here we define some automorphisms of the Yangian $Y(\mathfrak{gl}_{m|n})$ and homomorphisms between Yangians, so that we may refer to them in the next section.

Let 
$\omega_{m|n} : Y(\mathfrak{gl}_{m|n}) \to Y(\mathfrak{gl}_{m|n})$
be the automorphism defined by 
\begin{equation*} 
\omega : T(u) \mapsto T(-u)^{-1}.
\end{equation*}

Let $\tau: Y(\mathfrak{gl}_{m|n}) \to Y(\mathfrak{gl}_{m|n})$ be the anti-automorphism defined by 
\begin{equation*}
\tau(t_{ij}(u)) = t_{ji}(u) \times (-1)^{\pa{i}(\pa{j}+1)}.
\end{equation*}

Let $\rho_{m|n}: Y(\mathfrak{gl}_{m|n}) \to Y(\mathfrak{gl}_{n|m})$ be the isomorphism defined by
\begin{equation*}
\rho_{m|n} (t_{ij}(u)) = t_{m+n+1-i, m+n+1-j} (-u).
\end{equation*}

Let $\varphi_{m|n}: Y(\mathfrak{gl}_{m|n}) \hookrightarrow Y(\mathfrak{gl}_{m+k|n})$ be the inclusion which sends each generator $t_{ij}^{(r)} \in Y(\mathfrak{gl}_{m|n})$ 
to the generator $t_{k+i, k+j}^{(r)}$ in $Y(\mathfrak{gl}_{m+k|n})$.

Finally, let $\psi_k: Y(\mathfrak{gl}_{m|n}) \to Y(\mathfrak{gl}_{m+k|n})$ be the injective homomorphism defined by
\begin{equation}\label{psi}
\psi_k = \omega_{m+k|n} \circ \varphi_{m|n} \circ \omega_{m|n}.
\end{equation}
This last homomorphism is useful for studying quasideterminants so we discuss it in some detail with the following remarks.
 \begin{remark}\label{phil}
We can calculate $\psi_{k}(t_{ij}(u))$ explicitly for any $1 \le i,j \le m+n$ (see Lemma 4.2 of \cite{BK}) :
\begin{equation*}
\psi_{k} (t_{ij}(u)) = \left| \begin{array}{cccc} t_{11}(u) &\cdots &t_{1k}(u) &t_{1, k+j}(u)\\
\vdots &\ddots &\vdots &\vdots \\
t_{k1}(u) &\cdots &t_{kk}(u) &t_{k, k+j}(u)\\
t_{k+i, 1}(u) &\cdots &t_{k+i,k}(u) &\boxed{t_{k+i, k+j}(u)}
\end{array} \right|.
\end{equation*}
In particular, this means that 
for $k\ge 1$, we have
$\psi_{k} (d_1(u))= d_{k+1} (u),\;$ 
$\psi_{k} (e_1 (u))=e_{k+1}(u),$ and $\psi_{k}(f_1 (u))= f_{k+1}(u)$.

Furthermore, by (\ref{psi}), we have for any $k,l \ge 1$ that $\psi_k \circ \psi_l = \psi_{k+l}$, so we may generalise this observation to give for instance 
$\psi_{k}(d_l (u)) = d_{k+l}(u)$.
\end{remark}
\begin{remark}\label{phil2}
Notice that the map $\psi_k$ sends $t'^{\, (r)}_{ij} \in Y(\mathfrak{gl}_{m|n})$ to the element $t'^{\, (r)}_{k+i, k+j}$ in $Y(\mathfrak{gl}_{m+k|n})$. 
Thus the subalgebra $\psi_k (Y(\mathfrak{gl}_{m|n}))$  is generated by the elements $\{ t'^{\, (r)}_{k+s, k+t}\}_{s,t =1}^n$.  Then, by (\ref{useful}),
 all elements of this subalgebra commute with those of the subalgebra generated by the elements $\{ t_{ij}^{(r)} \}_{i,j =1}^{k}$.
 
By Remark \ref{phil}, this implies in particular that for any $i,j \ge 1$, the quasideterminants $d_i (u)$ and $d_j (v)$ commute.
\end{remark}
\section{The Quantum Berezinian}
The quantum Berezinian was defined by Nazarov \cite{Nazarov} and plays a similar role in the study of the Yangian $Y(\mathfrak{gl}_{m|n})$ as the quantum determinant does in
the case of the Yangian $Y(\mathfrak{gl}_{n})$ (see \cite{MNO}).
\begin{defn}\label{qbear}
The quantum Berezinian is the following power series with coefficients in the Yangian $Y(\mathfrak{gl}_{m|n})$:
{\small 
\begin{eqnarray*}
b_{m|n}(u)&:=& \nonumber \sum_{\tau \in S_m} \mathrm{sgn} (\tau) \,t_{\tau(1)1} (u) t_{\tau(2) 2} (u-1) \cdots t_{\tau(m)m} (u-m+1) \\
\nonumber&& \times \;\sum_{\sigma \in S_n}\! \mathrm{sgn}(\sigma) \, t'_{m+1, m+\sigma(1)} (u-m+1) 
 \cdots t'_{m+n, m+\sigma(n)} (u-m+n) \\
\end{eqnarray*}
}
\end{defn}
The first part of this expression for $b_{m|n}(u)$ is quite special and so is given its own notation:
\begin{equation*}
C_{m}(u) := \sum_{\tau \in S_m} \mathrm{sgn} (\tau) t_{\tau(1)1} (u) t_{\tau(2) 2} (u-1) \cdots t_{\tau(m)m} (u-m+1).
\end{equation*}
It is clear that $C_{m}(u)$ is an element of the subalgebra of $Y(\mathfrak{gl}_{m|n})$ generated by the set $\{t_{ij}^{(r)} \}_{1\le i,j \le m;r\ge 0}$.
This subalgebra is isomorphic to the Yangian $Y(\mathfrak{gl_m})$ associated to the Lie algebra $\mathfrak{gl}_m$ by the inclusion 
$Y(\mathfrak{gl}_m) \to Y(\mathfrak{gl}_{m|n})$ which send each generator $t_{ij}^{(r)}$ in $Y(\mathfrak{gl_m})$ to the generator of the same name in $Y(\mathfrak{gl}_{m|n})$.  
Moreover, $C_{m}(u)$ is in fact the image under this map of the \emph{quantum determinant} of the smaller Yangian $Y(\mathfrak{gl_m})$ (see \cite{BK},
\cite{MNO}).
Then it is well known (see Theorem 2.32 in \cite{Molev}) that we have the alternative expression:
\begin{equation*}
C_{m}(u) = d_{1}(u) d_{2}(u-1) \cdots d_{m}(u-m+1).
\end{equation*}
We can extend this observation as follows:
\begin{theorem}We have the following alternative expression for the quantum Berezinian:
\begin{eqnarray} 
b_{m|n}(u)
&=& \nonumber d_{1}(u)\, d_{2}(u-1) \cdots d_{m}(u-m+1) \\
\nonumber&& \times\, d_{m+1}(u-m+1)^{-1} \cdots d_{m+n} (u-m+n)^{-1}.
\end{eqnarray}
\end{theorem}
\begin{proof}
Notice that the second part of the expression for $b_{m|n}(u)$ in Definition~\ref{qbear} is the image under the isomorphism 
$\rho_{n|m} \circ \omega_{n|m}: Y(\mathfrak{gl}_{n|m}) \to Y(\mathfrak{gl}_{m|n})$ of 
{\small
\begin{equation}\label{bear}
\sum_{\sigma \in S_n}\! \mathrm{sgn}(\sigma) \, t_{n, \sigma(n)} (u-m+1) \cdots t_{2, \sigma(2)}(u-m+n-1) \, t_{1, \sigma(1)} (u+m-n)
\end{equation}
}where in this expression (\ref{bear}) we are following the usual convention for denoting generators \emph{in the Yangian $Y(\mathfrak{gl}_{n|m})$}.  
We recognise (by comparing with (8.3) of \cite{BK} for example) that the expression (\ref{bear}) is in fact $C_{n}(u-m+n)$, the image of the quantum determinant of
$Y(\mathfrak{gl}_n)$ under the natural inclusion $Y(\mathfrak{gl}_n)~\hookrightarrow~Y(\mathfrak{gl}_{n|m})$. 
So in order to verify the claim we must calculate the image of $C_{n}(u-m+n)$ under this map explicitly in terms of our quasideterminants $d_i(v)$.
Applying Proposition 1.6 of \cite{GR}, we find that the image of $d_i(v)$~in~$Y(\mathfrak{gl}_{n|m})$ is $(d_{m+n+1-i}(v))^{-1}$~in~$Y(\mathfrak{gl}_{m|n})$. This gives the
desired result. 
\end{proof}  
The following theorem is a result of Nazarov \cite{Nazarov}.  We give a new proof.
\begin{theorem}
The coefficients of the quantum Berezinian (\ref{qbear}) are central in the algebra $Y(\mathfrak{gl}_{m|n})$.
\end{theorem}
\begin{proof}
By Remark \ref{phil2}, we already know that the quantum Berezinian $b_{m|n}(u)$ commutes with $d_{i}(v)$ for $1 \le i \le m+n$.  
So our problem reduces to showing that $b_{m|n}(u)$ commutes with $e_{i}(v)$ and with $f_{i}(v)$ for each $i$ between $1$ and $m+n-1$.  We proceed by breaking this problem into three cases.

\emph{Case 1: $1\le i \le m-1$.} In this case, $e_{i}(v)$ commutes with $C_{m}(u)=d_{1}(u)~\cdots~d_{m}(u-m+1)$ by Theorem 7.2 in \cite{BK}.  On the
other hand, $e_{i}(v)$ is an element of the subalgebra generated by $\{t_{jk}^{(r)} \}_{1 \le j,k \le m}$ and thus by Remark \ref{phil2} commutes with  
$d_{m+s}(u-m+s)^{-1}=t'_{m+s, m+s}(u-m+s)$ for $1\le s\le n$. Now we may use the anti-automorphism $\tau$ to show that $f_{i}(v)$ also commutes with the quantum Berezinian in this case
because $\tau (e_{i}(v)) = f_{i}(v)$ for $1\le i \le m-1$.

\emph{Case 2: $m+1\le i \le m+n-1$.}
Applying Propositions 1.6 and 1.4 of \cite{GR} in turn to $f_{i}(v)$, we find an alternative expression:
{\small \begin{eqnarray*}
f_{i}(v) &=& -  \left| \begin{array}{cccc} \boxed{t'_{i+1,i}(v)} & t'_{i+1, i+2}(v)& \cdots & t'_{i+1, m+n}(v)\\
t'_{i+2,i}(v) & t'_{i+2, i+2}(v)&\cdots & \\
\vdots  &\vdots &\ddots &\vdots\\
t'_{m+n,i}(v) &t'_{m+n, i+2}(v)& \cdots & t'_{m+n, m+n}(v)\end{array} \right|\, \left| \begin{array}{ccc} \boxed{t'_{i+1, i+1}(v)}&\cdots& t'_{i+1, m+n}(v)\\
\vdots  & &\vdots\\
t'_{m+n, i+1}(v) &\cdots& t'_{m+n, m+n}(v)
\end{array} \right|^{-1}.
\end{eqnarray*} 
}
Thus, we find that for $m+1\le i \le m+n-1$,
\begin{equation*}
\rho_{n|m} \circ \omega_{n|m} (-f_{m+n-i}(v)) = e_{i}(v),
\end{equation*} 
and similarly
\[
\rho_{n|m} \circ \omega_{n|m} (-e_{m+n-i}(v)) = f_{i}(v).
\]
We apply this isomorphism to the results of Case 1 in the Yangian $Y(\mathfrak{gl}_{n|m})$. 
This shows that $e_{i}(v)$ and $f_{i}(v)$ commute with the quantum Berezinian in the case where $m+1 \le i \le m+n-1$.  

\emph{Case 3: i=m.} We begin by considering the Yangian $Y(\mathfrak{gl}_{1|1})$.  For this algebra we have that the quantum Berezinian is $b_{1|1} (u) = d_{1}(u) d_{2}(u)^{-1}$ and we would like to show
that it commutes with $e_{1}(v)$ and $f_{1}(v)$.   So it will suffice to show
\begin{equation}\label{star}
d_{1}(u) e_{1}(v) d_{2}(u) = d_{2} (u) e_{1} (v) d_{1}(u).
\end{equation}
We have
{\small  \begin{eqnarray}\label{matrix}
\left( \begin{array}{cc} t_{11} (u)\! & t_{12}(u)\\
t_{21} (u)\! & t_{22}(u)
\end{array} \right)  &=& \left( \begin{array}{ll} d_{1}(u) & d_{1}(u)\,e_{1}(u)\\
f_{1}(u)d_{1}(u) & f_1 (u)d_1 (u) e_1 (u) + d_2(u)
\end{array} \right) \\
\left( \begin{array}{ll} t'_{11} (v)\! & t'_{12}(v)\\
t'_{21} (v)\! & t'_{22}(v)
\end{array} \right) \label{matrixinverse}
&=& \left( \begin{array}{ll} d_{1}(v)^{-1}\! + e_{1}(v) d_{2}(v)^{-1} f_{1}(v) & -e_{1}(v)\, d_{2} (v)^{-1}\!\!\!\\
- d_{2}(v)^{-1} f_{1}(v) & \phantom{-}d_2(v)^{-1}
\end{array} \right).
\end{eqnarray} 
}
An application of (\ref{useful}) gives 
\begin{equation*}
(u-v) [t_{11} (u), t'_{12}(v)] \;=\; t_{11}(u) t'_{12}(v) + t_{12}(u) t'_{22} (v).
\end{equation*}
Substituting in the expressions from (\ref{matrix}) and (\ref{matrixinverse}) then cancelling $d_{2} (v)$, this gives
\begin{equation*}
(u-v) [ d_{1} (u), e_{1} (v)] \; =\; d_{1} (u) (e_{1} (v) -e_{1}(u)).
\end{equation*}
Similarly, by considering the commutator $[t_{12}(u), t_{22}'(v)]$, we derive the relation
\begin{equation*}
(u-v) [ d_{2} (u), e_{1} (v)] \; =\; d_{2} (u) (e_{1} (v) -e_{1}(u)).
\end{equation*}
We rewrite these relations to find
\begin{eqnarray*}
(u-v) e_{1} (v) d_{1} (u) &=& (u-v-1) d_{1}(u) e_{1} (v) + d_{1} (u) e_{1} (u),\\
(u-v) e_{1} (v) d_{2} (u) &=& (u-v-1) d_{2}(u) e_{1} (v) + d_{2} (u) e_{1} (u),
\end{eqnarray*}
and by considering these expressions we see that (\ref{star}) holds.  Similarly, the quantum Berezinian commutes with $f_{1}(v)$.

Now we return our attention to the general Yangian $Y(\mathfrak{gl}_{m|n})$.  By similar appeals to Remark~\ref{phil2} as in the first case, we see that $e_{m}(v)$ and 
$f_{m}(v)$ commute with $d_{1} (u) \cdots d_{m-1}(u-m+2)$ and with $\;d_{m+2}(u-m+2)^{-1} \cdots~d_{m+n}(u-~m~+~n)^{-1}$.  
So we need only show they commute with $d_{m}(u-m+1)d_{m+1}(u-m+1)^{-1}$.  
This follows immediately when we apply the homomorphism $\psi_{m-1}$  to the identity (\ref{star}) in $Y(\mathfrak{gl}_{1|1})$.
\end{proof}
\bibliographystyle{plain}
\def\cprime{$'$} \def\cprime{$'$} \def\cprime{$'$} \def\cprime{$'$}
  \def\cprime{$'$}

\end{document}